\documentclass[12pt,reqno]{amsart}
\usepackage{fullpage}
\usepackage{amsfonts}
\usepackage{amssymb}
\usepackage{times}
\usepackage{graphicx}
\usepackage{float} 
\begin{document}
\title{Iterative Aggregation Method for Solving Principal Component Analysis Problems}
\maketitle
\begin{minipage}[H]{0.6\textwidth}
{\bf  Vitaly Bulgakov}
\end{minipage}
\begin{minipage}[H]{0.6\textwidth}
BULGAKOV\_V@YAHOO.COM
\end{minipage}

\renewcommand\abstractname{\textbf{Abstract}}

\begin{abstract}

Motivated by the previously developed multilevel aggregation method for solving structural analysis problems a novel two-level aggregation approach for efficient iterative solution of Principal Component Analysis (PCA) problems is proposed. The course aggregation model of the original covariance matrix is used in the iterative solution of the eigenvalue problem by a power iterations method. The method is tested on several data sets consisting of large number of text documents.

\end{abstract}

\section{\bf Introduction}
\label{first section}

This work was envisioned as application of the {\it multilevel aggregation method} \cite{ref-mlam} developed by the author back in 90s to PCA problems.
Multilevel aggregation method was an extension of well-known {\it multigrid} methods\cite{ref-mgm} from boundary value problems to general structural analysis problems which brought it to the class of {\it algebraic multigrid} methods. The idea of the aggregation method was to use some naturally constructed course model of the original finite element approximation of a structure which provides a fast convergence for iterative methods for solving large algebraic systems of equations. One of applications of this method was an iterative solution of large eigenvalue problems arising in structural natural vibration and buckling analyses \cite{ref-vib-buc}. In these problems a sought set of lowest vibration modes can be thought of as {\it principal components} of structure behavior. An obvious similarity with PCA was a turning point to start looking for a proper way to create an aggregation model for data matrix approximation and use it for efficient solution of PCA problems.     

In this study PCA\cite{ref-pca} is applied to and the method is tested on text analysis problems. A tested data set consists of documents each of which produces an N-dimensional vector stored as a column of a data matrix which values are term frequencies. Our raw data comes in the form of text files from data sets such as medical abstracts and news groups. The purpose of PCA is to iteratively compute a set of highest eigenvalues and corresponding eigenvectors of the covariance matrix. Covariance matrix is never formed explicitly. The main operation is multiplication of large sparse data matrix or its transpose by a vector. The course aggregation model of the original covariance matrix is used in the iterative solution of the eigenvalue problem. Original covariance matrix and its approximation of small size assumes similarity of leading eigenvalues and eigenvectors. This fact allows fast convergence of subspace iterations at minimal additional computational cost.  

For numerical experiments we use R language which is rich of linear algebra, statistical and graphical packages. 

\section{\bf PCA problem formulation}
PCA in multivariate statistics is widely used as an effective way to perform unsupervised dimension reduction. The essence of this method lies in using Singular Value Decomposition (SVD) which provides the best low rank approximation to original data according to Eckart-Young theorem \cite{ref-ey}. 
Let $n$ data points in $m$ dimensional space be contained in the data matrix which is assumed already centered  around the origin for computational stability
\begin{equation}
\label{1}
(x_1, x_2,...,x_n) = X
\end{equation}
Then covariance matrix is 

\begin{equation}
\label{2}
A = XX^T
\end{equation}

Let ($\lambda_k$, $\phi_k$) be an eigenpair of $A$, where eigenvectors $\phi_k$ define principal directions.

\section{\bf Aggregation model}
In order to create an aggregation model we divide the entire set of data vectors $x_i$ into $n_0$ clusters using some similarity criteria where $n_0 << n$. We will explain later how we do clustering. We assume that all vectors within the cluster are similar and a single representative of a cluster is an average of all vectors $x_i$ where $i \in cluster_k$ or for cluster $k$ we have
\begin{equation}
\label{3}
x_k^0 = 1/{dim_k}*\sum_{i \in cluster_k}x_i
\end{equation}
Transformation of matrix $X$ to $X_0$ is done using matrix $R$ which we call {\it aggregator}
\begin{equation}
\label{4}
X_0 = XR 
\end{equation} 
where $R[i,k] =$ if  $i \in cluster_k$  then  $1/{dim_k}$ else $0$. $X_0$ is of size $(m,n_0)$. Approximation $A_0$ of covariance matrix $A$ is 
\begin{equation}
\label{5}
A_0 = X_0 X_0^T = XRR^TX^T
\end{equation}
 Formally matrix $A_0$ is of the same size as $A$ but has a much lower rank. We do not need to use form \eqref {5} for computations. For matrix vector multiplication we rather use sparse   
matrix $X_0$ which according to \eqref{3} is constructed by simple averaging of vectors inside a cluster and 
\begin{equation}
\label{6}
A_0 v = X_0 X_0^T v
\end{equation}
Therefore $A_0 v$ requires $O(mn_0)$ operations which is much lower than $O(mn)$ operations required for $Av$. We also expect and this is confirmed by numerical experiments that convergence of iterative methods for solving partial eigenvalue problem for $A_0$ is faster than that for $A$.

There are quite a few clustering techniques known as computationally efficient. Besides since we need clustering as an auxiliary procedure we do not need highly accurate clustering results. In this study we use K-means clustering algorithm \cite{ref-kmeans} which became very popular in data mining, unsupervised classification, etc. and which converges quickly to a local optimum. Our experience says that the aggregated problem with a small number of clusters provides a good resemblance of the original and approximated covariance matrices in terms of first (highest) eigenvalues which is important for the iterative method described below. In Figure 1 this resemblance is demonstrated where we show distribution of first 10 eigenvalues of both matrices where the data matrix $X$ was obtained by processing 2014 documents of "Cardiovascular Diseases Abstracts" corpus. Matrix $X_0$ was obtained by K-means method with 10 clusters.

\begin{figure}[H]
  \includegraphics[width=0.5\textwidth]{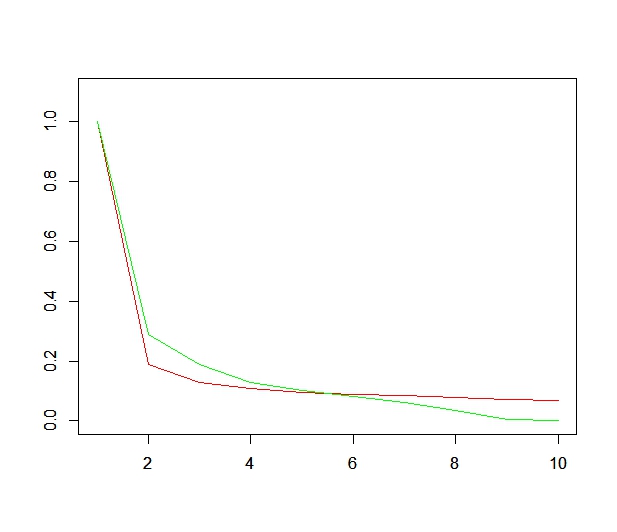}
  \caption{Distribution of 10 first eigenvalues of the original and approximate covariance matrices for 16058 by 2014 data matrix and 10 clusters}
  \label{fig:1}
\end{figure}

\section{\bf Iterative method}
We use power iteration method \cite{ref-pi} for for solving auxiliary aggregated eigenvalue problem and a modified power method for solving the original eigenvalue problem. This method is also known as subspace iteration when used to simultaneously iterate a set of eigenvectors. One iteration of the power algorithm consists of the following steps:
\begin{equation}
\label{7}
\begin{aligned}
for \; i \; =  \; 1 \; to \; l \;:\;\;\;
\tilde{u_i}^{k+1} = \frac{1}{\|Au_i^k\|} * Au_i^k \\
u_1^{k+1},...,u_l^{k+1} = orthonorm (\tilde{u_1}^{k+1},...,\tilde{u_l}^{k+1}) \\
with \; approximation \; of \; eigenvalues \\
\lambda_i^k = \frac{(Au_i^k,u_i^k)}{(u_i^k,u_i^k)}
\end{aligned}
\end{equation}
which starts with a set of $l$ initial approximations of first eigenvectors $(u_1^0,u_2^0,...,u_l^0) = U^0$. The key property of the power method is that if approximation $u_i^0$ is spanned by matrix $A$ eigenvectors subspace, then after $k$ multiplications of matrix $A$ by this vector the linear combination of eigenvectors will be weighted by $\lambda_i$ to the power $k$ which gives boost to terms corresponding to highest eigenvalues:

 \begin{equation}
\label{7.1}
\begin{aligned}
A^ku = \sum_k c_i\lambda_i^k\phi_i
\end{aligned}
\end{equation}  

In the method proposed for the first $l$ principal directions of PCA we will need first $k$ orthonormal eigenvectors of $A_0$ $q_1, q_2, ..., q_k$ where $k >= l$. These vectors can be obtained by algorithms \eqref{7}. We will also need matrix $P_i$
\begin{equation}
\label{8}
P_i = q_i  q_i^T 
\end{equation} 
Since $q_i^Tq_j = \delta_{i,j}$ and $P_i  P_i = P_i$, it is a projector to the subspace of $i$-th eigenvector of $A_0$. We will modify method \eqref {7} using this projector in the following manner:

 \begin{equation}
\label{9}
\begin{aligned}
for \; i \; =  \; 1 \; to \; l \;:\;\;\;
\tilde{u}_i^{k+1} = \frac{1}{\|Bu_i^k\|} * Bu_i^k \;\; where \\
Bu_i^k = Au_i^k + \alpha_i * P_iAP_iu_i^k  \;\; and \;\; \alpha_i => min\|A\tilde{u}_i^{k+1} - \lambda_i^k\tilde{u}_i^{k+1}\| \\
u_1^{k+1},...,u_l^{k+1} = orthonorm (\tilde{u_1}^{k+1},...,\tilde{u_l}^{k+1}) \\
with \; approximation \; of \; eigenvalues \\
\lambda_i^k = \frac{(Au_i^k,u_i^k)}{(u_i^k,u_i^k)}
\end{aligned}
\end{equation}

This approach can be thought of as "help" to the power iteration method to converge on the subspace of eigenvectors of the aggregated problem. The intuition for that is similarity of first eigenvectors and eigenvalues of the original and aggregated problem if clustering is done properly. Let $u = \sum c_i\phi_i$ where $\phi_i$ are eigenvectors of the original covariance matrix $A$ and $P_k$ be a projector on subspace of $\phi_k$. Then

 \begin{equation}
\label{10}
\begin{aligned}
Au + \alpha P_kAP_k u= \sum_{i \neq k} c_i\lambda_i\phi_i + c_k\lambda_k(1+\alpha)\phi_k
\end{aligned}
\end{equation}

If $\alpha$ is chosen big then the second term of this expression dominates over the first term thus providing convergence for $\phi_k$ in one iteration step. $\alpha$ can be derived from the condition stated in \eqref {9}:

\begin{equation}
\label{11}
\begin{aligned}
\Phi(\tilde{u}_i^{k+1}) =  \|A\tilde{u}_i^{k+1} - \lambda_i^k\tilde{u}_i^{k+1}\| \\
\alpha => min \Phi => \frac{d \Phi}{d \alpha} = 0 
\end{aligned}
\end{equation}

This equation leads to the quadratic equation for $\alpha$. Omitting indexes and skipping details we arrive at the following expression for $\alpha$

\begin{equation}
\label{12}
\begin{aligned}
\alpha = -\frac{(A^2u,AFu) - 2 \lambda (A^2u,Fu) + \lambda^2(Au,Fu)}{(Au,AFu) - 2 \lambda (Au,Fu) + \lambda^2(Fu,Fu)}
\end{aligned}
\end{equation}  

where $F = PAP$.  

We note that as you can see from \eqref {11} $\alpha$ is chosen from the previous step to simplify computations. This can also be justified by the fact that eigenvalues converge faster than eigenvectors. Detailed algorithm discussion is out of scope of this paper. We just mention here that all operations with matrix $A$ are reduced to the matrix vector multiplications of the sparse data matrix $X$ or its transpose $X^T$. 

\section{\bf Numerical experiments}

For numerical experiments we used two data sets. The fist one is "Cardiovascular Diseases Abstracts" which is a set where each abstract is an individual document. The data matrix $X$ size is 16058 by 2014 where the first value is the total number of terms and the second one is the number of documents. We searched for 10 first eigenvalues of the covariance matrix $A=XX^T$ and used 10 clusters for constructing auxiliary aggregation problem $A_0 = X_0X_0^T$. So the size of this problem is more than 201 times lower than that for the original problem.

The problem is solved using algorithm \eqref {9}.  Figure 2 shows changes of parameter $\alpha$ for the first three eigenvectors. As expected the biggest contribution of projectors \eqref {8} is observed in first iterations to suppress errors caused by initial eigenvector guesses. After some number of iteration contribution of projectors is getting smaller while eigenvectors are getting more accurate. 

\begin{figure}[H]
  \includegraphics[width=0.5\textwidth]{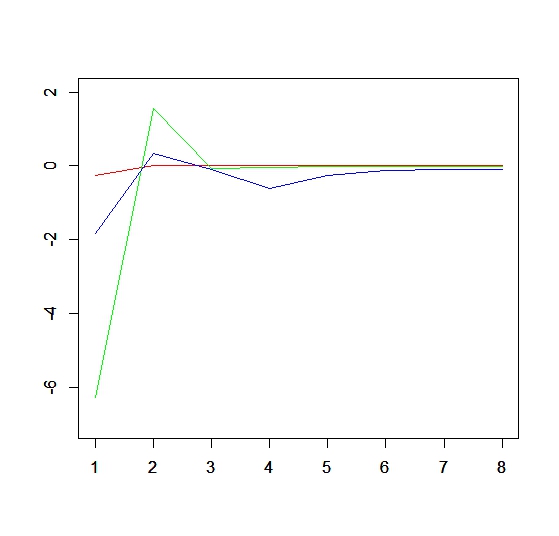}
  \caption{$\alpha$ changes with iteration number. $\alpha_1$ (for eigenvector 1) - red, $\alpha_2$ (for eigenvector 2) - green, $\alpha_3$ (for eigenvector 3) - blue}
  \label{fig:2}
\end{figure}

We measure convergence of eigenvalues through $Error_1 = \|\Lambda^{k+1} - \Lambda^k\|_F/\|\Lambda^k\|_F$ and convergence of eigenvectors by the residual matrix through    
$Error_2 = \|AU^k - U^{k}\Lambda^k\|_F$ 
where $\|\|_F$ is a matrix Frobenius norm,  $U^k$ consists of orthonormal vectors $u_1^{k},...,u_l^{k}$ which are approximations of the eigenvectors and 
$\Lambda^k$ is a diagonal matrix of approximations of eigenvalues. Errors graph is demonstrated in Figure 3.

A good convergence rate of the iterative process is demonstrated. After 40 iterations we got $Error_1 = 0.00038$ and $Error_2 = 0.0017$. 

\begin{figure}[H]
  \includegraphics[width=0.5\textwidth]{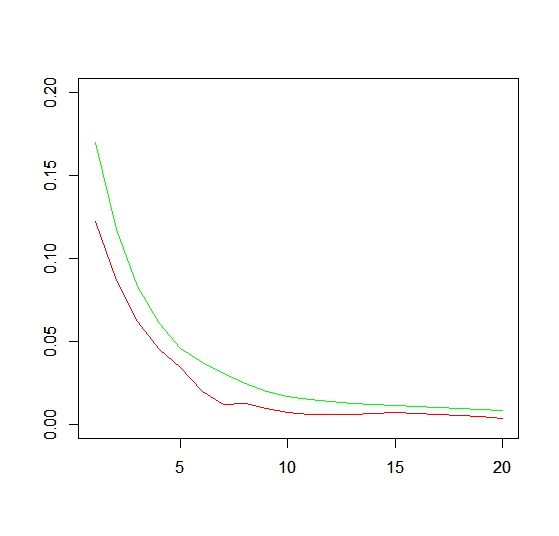}
  \caption{Eigenvalues and Eigenvectors convergence for "Cardiovascular Diseases Abstracts" data set. $Error_1$ - red, $Error_2$ - green.}
  \label{fig:3}
\end{figure}

The second corpus was ''talk politics'' set from the news groups. Size of this problem is 13511 (terms) by 1171 (documents).  We searched for 10 first eigenvalues of the covariance matrix and used 10 clusters again. The quality of the clustering aggregated model can be viewed by comparing eigenvalues of the original and aggregated covariance matrices. Figure 4 demonstrates a good resemblance of eigenvalues distribution. 
Convergence graph is demonstrated in Figure 5. After 40 iterations we got $Error_1 = 0.00044$ and $Error_2 = 0.00049$.

\begin{figure}[H]
  \includegraphics[width=0.5\textwidth]{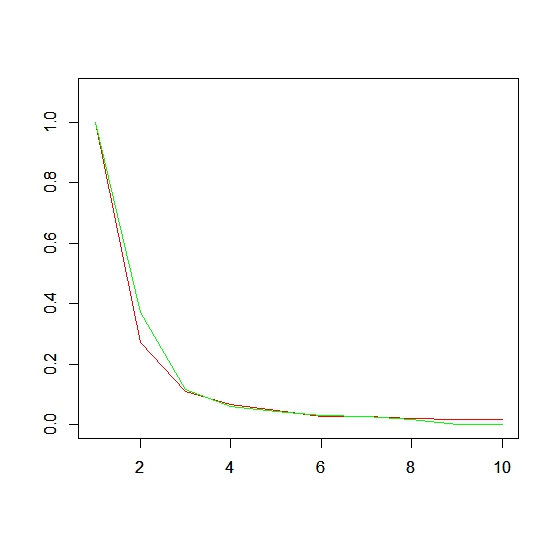}
  \caption{Distribution of 10 first eigenvalues of the original and approximate covariance matrices for 13511 by 1171 data matrix and 10 clusters}
  \label{fig:4}
\end{figure}

\begin{figure}[H]
  \includegraphics[width=0.5\textwidth]{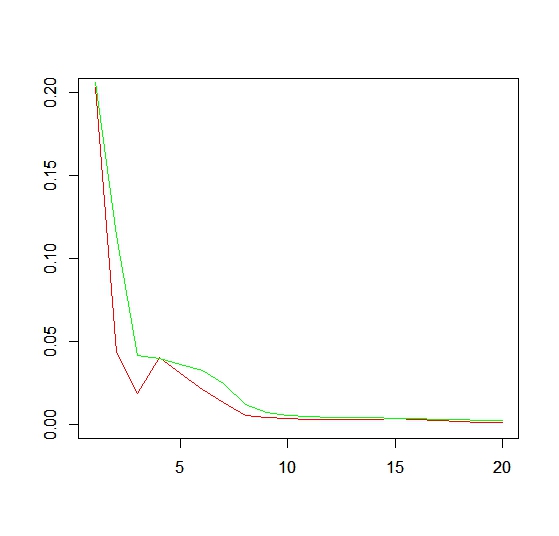}
  \caption{Eigenvalues and Eigenvectors convergence for ''News Group'' corpus. $Error_1$ - red, $Error_2$ - green.}
  \label{fig:5}
\end{figure}

\end{document}